\documentclass[reqno, 11pt]{article}

\usepackage{enumerate}
\usepackage{needspace}
\usepackage{latexsym}
\usepackage{amsmath}
\usepackage{amsfonts}
\usepackage{amsthm}
\usepackage{amssymb,mathtools}
\usepackage{newlfont}
\usepackage{graphics}
\usepackage{color}
\usepackage{float}
\usepackage{diagbox}
\usepackage{tocloft}
\usepackage{titlesec}
\usepackage{booktabs,array,longtable,enumitem,microtype}
\usepackage{tikz}
\usepackage{xcolor}
\usepackage{extpfeil}
\usepackage{centernot}
\usepackage{titlesec}
\usepackage{graphicx}
\usepackage{subfig}
\usetikzlibrary{arrows.meta, calc}

\usepackage[linesnumbered,ruled,vlined]{algorithm2e}
\usepackage{url}

\usepackage[margin=1in]{geometry}

\usepackage[T1]{fontenc}
\usepackage{lmodern}
\usepackage{rotating}
\usepackage{multirow}
\usepackage{extarrows}
\usepackage[sort,compress,numbers]{natbib}
\usepackage[utf8]{inputenc}
\usepackage{listings}
\usepackage{aliascnt}
\usepackage[pagebackref,colorlinks=true,linkcolor=blue,citecolor=red,urlcolor=blue]{hyperref}
\usepackage[nameinlink,noabbrev,capitalize]{cleveref}
\numberwithin{equation}{section}

\usepackage{enumitem}
\setlist[enumerate]{label=\textup{(\roman*)},leftmargin=*}
\newtheorem{theorem}{Theorem}[section]
\newaliascnt{proposition}{theorem}
\newtheorem{proposition}[proposition]{Proposition}
\aliascntresetthe{proposition}
\newaliascnt{lemma}{theorem}
\newtheorem{lemma}[lemma]{Lemma}
\aliascntresetthe{lemma}
\newaliascnt{corollary}{theorem}

\aliascntresetthe{corollary}

\theoremstyle{definition}
\newaliascnt{definition}{theorem}

\aliascntresetthe{definition}
\newaliascnt{remark}{theorem}

\aliascntresetthe{remark}
\newaliascnt{example}{theorem}

\aliascntresetthe{example}
\newaliascnt{conjecture}{theorem}
\newtheorem{conjecture}[conjecture]{Conjecture}
\aliascntresetthe{conjecture}
\newaliascnt{question}{theorem}

\aliascntresetthe{question}
\newaliascnt{problem}{theorem}

\aliascntresetthe{problem}
\newaliascnt{alg}{theorem}

\aliascntresetthe{alg}

\allowdisplaybreaks[4]

\SetKwInput{KwInput}{Input}                
\SetKwInput{KwOutput}{Output}              

\DeclareMathOperator{\vol}{vol}

\newcommand{\dd}{\,\mathrm d}

\newcommand{\ee}{\mathrm e}
\newcommand{\Q}{\mathbb Q}
\newcommand{\C}{\mathbb C}

\newcommand{\Ztwo}{\mathrm{Z2}}
\newcommand{\Zthree}{\mathrm{Z3}}

\title{Two families of Young tableaux with periodic walls:\\ Counting formulas and differential transcendence}

\author{Feihu Liu$^{\color{blue} \dag}$, Jiahang Liu$^{\textcolor{blue} \ddag}$\thanks{Corresponding author}, and Guoce Xin$^{\color{blue} \S}$
\\[2mm]
{\small $^{\color{blue} \dag}$ Center for Combinatorics, LPMC,}\\[-0.8ex]
{\small Nankai University, Tianjin 300071, P.R.~China}\\
{\small $^{\color{blue} \ddag, \S}$ School of Mathematical Sciences,}\\[-0.8ex]
{\small Capital Normal University, Beijing 100048, P.R.~China}\\
{\small {\color{blue} $^\dag$} Email address: liufeihu7476@163.com}\\
{\small {\color{blue} $^\ddag$} Email address: jiahang-liu@foxmail.com}\\
{\small {\color{blue} $^\S$} Email address: guoce\_xin@163.com}
}

\date{\today}

\begin{document}

\maketitle

\begin{abstract}
Banderier, Marchal, and Wallner studied the enumeration problem for Young tableaux with walls. These are similar to standard Young tableaux, except that local decreases are allowed at some walls. In 2021, Banderier and Wallner counted Young tableaux with periodic walls, in which the elementary blocks are Young diagrams of shape $(2,2)$. They obtained exact counting formulas for $30$ of the $32$ classes up to symmetry. In this paper, we first derive recurrences for counting the remaining two classes of Young tableaux with periodic walls. We then prove that their counting exponential generating functions are differentially transcendental. Consequently, they are neither D-algebraic nor D-finite. This disproves a conjecture of Banderier and Wallner.
\end{abstract}

\noindent
\begin{small}
\emph{2020 Mathematics subject classification}: Primary 05A15; Secondary 06A07, 34M15.
\end{small}

\noindent
\begin{small}
\emph{Keywords}: Young tableaux with walls; Differentially transcendental; D-algebraic; D-finite; Linear extensions; Exponential generating function.
\end{small}


\section{Introduction}

Banderier, Marchal, and  Wallner~\cite{C.Banderier18,banderier2021young} considered a variation of standard Young tableaux, called \emph{Young tableaux with walls}. It is very difficult to find a counting formula for such Young tableaux.
We now introduce the necessary definitions.

Let $\lambda=(\lambda_1,\lambda_2,\ldots)$ be a partition of $m$, denoted $\lambda \vdash m$.
The \emph{Young diagram} of $\lambda$ is obtained by drawing a collection of left-justified arrays of $m$ cells, where row $i$ has $\lambda_i$ cells reading from bottom to top. We also denoted the Young diagram of $\lambda$ by $\lambda$.
A \emph{standard Young tableau} of $\lambda$ is an assignment of $1,2,\ldots,m$ to the cells of the Young diagram, each number appearing once, such that rows are increasing from left to right, and columns are increasing from bottom to top.
We refer to \cite[Section 1.7]{StanleyEC1} for details of these definitions.

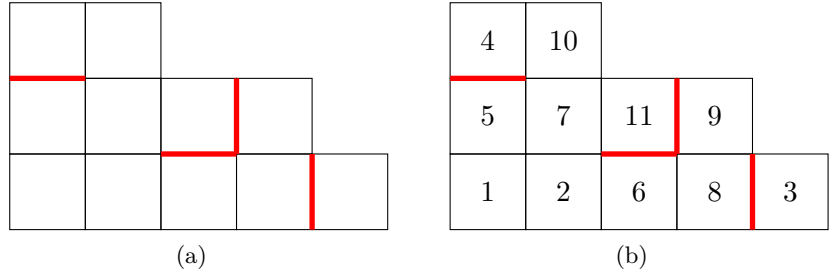
\begin{figure}[htbp]
\centering
\subfloat[]{
\begin{tikzpicture}[
    cell/.style={draw, minimum size=1cm, anchor=center},
    wall/.style={line width=2pt, red}]
\node[cell] (a1) at (0,0) {};
\node[cell] (a2) at (1,0) {};

\node[cell] (b1) at (0,-1) {};
\node[cell] (b2) at (1,-1) {};
\node[cell] (b3) at (2,-1) {};
\node[cell] (b4) at (3,-1) {};

\node[cell] (c1) at (0,-2) {};
\node[cell] (c2) at (1,-2) {};
\node[cell] (c3) at (2,-2) {};
\node[cell] (c4) at (3,-2) {};
\node[cell] (c5) at (4,-2) {};

\draw[wall] (-0.5,-0.5) -- (0.5,-0.5);
\draw[wall] (2.5,-0.5) -- (2.5,-1.5);

\draw[wall] (1.5,-1.5) -- (2.5,-1.5);
\draw[wall] (3.5,-1.5) -- (3.5,-2.5);
\end{tikzpicture}}
\hspace{15pt}
\subfloat[]{
\begin{tikzpicture}[
    cell/.style={draw, minimum size=1cm, anchor=center},
    wall/.style={line width=2pt, red}]
\node[cell] (a1) at (0,0) {4};
\node[cell] (a2) at (1,0) {10};

\node[cell] (b1) at (0,-1) {5};
\node[cell] (b2) at (1,-1) {7};
\node[cell] (b3) at (2,-1) {11};
\node[cell] (b4) at (3,-1) {9};

\node[cell] (c1) at (0,-2) {1};
\node[cell] (c2) at (1,-2) {2};
\node[cell] (c3) at (2,-2) {6};
\node[cell] (c4) at (3,-2) {8};
\node[cell] (c5) at (4,-2) {3};

\draw[wall] (-0.5,-0.5) -- (0.5,-0.5);
\draw[wall] (2.5,-0.5) -- (2.5,-1.5);

\draw[wall] (1.5,-1.5) -- (2.5,-1.5);
\draw[wall] (3.5,-1.5) -- (3.5,-2.5);
\end{tikzpicture}}
\caption{A Young building and a Young tableau with walls.}\label{tab:youngExample}
\end{figure}

A \emph{Young building} $\mathcal{B}$ of shape $\lambda$ is a pair $\mathcal{B}=(\lambda,\mathcal{W})$, where $\mathcal{W}$ is a subset of the edges of Young diagram of $\lambda$. The elements of $\mathcal{W}$ are called \emph{walls}.
A \emph{Young tableau with walls over $\mathcal{B}$} is a filling of the cells of $\lambda$
by labels $1,2,\dots, m=|\lambda|$, such that each label appears exactly once and the labels are increasing along each row and column (conditions for standard Young tableaux), but two labels separated by a wall in $\mathcal{W}$ need not be increasing.
Denote by $\mathrm{YT}(\lambda,\mathcal{W})$ (or $\mathrm{YT}(\mathcal{B})$) the set of all such tableaux.
It reduces to the set of standard Young tableaux of shape $\lambda$ when $\mathcal{W}$ is empty.
See \cref{tab:youngExample} for an example of Young tableaux with walls of shape $\lambda=(5,4,2)$ in $\mathrm{YT}(\lambda,\mathcal{W})$. The walls in $\mathcal{W}$ are depicted as bold red edges.

Banderier, Marchal, and Wallner~\cite{C.Banderier18} studied the enumeration of certain rectangular Young tableaux with walls using density methods and bijections.
Through the enumeration of two classes of Young tableaux with walls, Lin et al. \cite{LinEtAl2026} proved the well-known tree-child network conjecture of Pons and Batle \cite{PonsBatle2021}.
The solution of this conjecture was facilitated by the connection established by Chang et al. \cite[Remark 3.3]{CFLWY} between tree-child networks and Young tableaux with walls.

Banderier and Wallner~\cite{banderier2021young} considered Young tableaux with periodic walls and obtained several counting formulas. By employing the kernel method and the lattice path enumeration theory, Liu and Xin~\cite{LiuXin-YTW} obtained the generating function for a class of two-row Young tableaux with periodic walls.

Let $\mathcal{B}=(\lambda,\mathcal{W})$ be a building of rectangular shape.
Define $\mathcal{B}^n$ to be the Young building formed by connecting $n$ blocks of $\mathcal{B}$ horizontally.
Then Young tableaux with walls over $\mathcal{B}^n$ are called \emph{Young tableaux with periodic walls}.
Let $\mathrm{YT}(\mathcal{B}^n)$ denote the set of all such tableaux.
Note that $\mathcal{W}$ is allowed to contain the rightmost edges of $\lambda$.
\cref{tab:youngZ2} illustrates a Young building $\mathcal{B}$ on the left and a Young building $\mathcal{B}^3$ with periodic walls on the right.

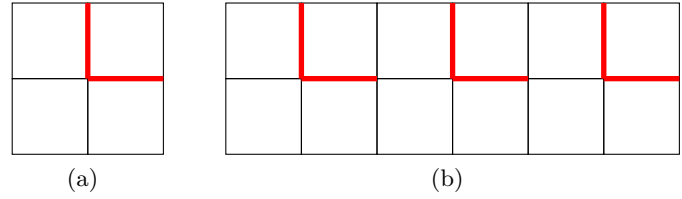
\begin{figure}[htbp]
\centering
\subfloat[]{
\begin{tikzpicture}[
    cell/.style={draw, minimum size=1cm, anchor=center},
    wall/.style={line width=2pt, red}]
\node[cell] (a1) at (0,0) {};
\node[cell] (a2) at (1,0) {};

\node[cell] (b1) at (0,-1) {};
\node[cell] (b2) at (1,-1) {};

\draw[wall] (0.5,0.5) -- (0.5,-0.5);
\draw[wall] (0.5,-0.5) -- (1.5,-0.5);
\end{tikzpicture}}
\hspace{15pt}
\subfloat[]{
\begin{tikzpicture}[
    cell/.style={draw, minimum size=1cm, anchor=center},
    wall/.style={line width=2pt, red}]
\node[cell] (a1) at (0,0) {};
\node[cell] (a2) at (1,0) {};
\node[cell] (a3) at (2,0) {};
\node[cell] (a4) at (3,0) {};
\node[cell] (a5) at (4,0) {};
\node[cell] (a6) at (5,0) {};

\node[cell] (b1) at (0,-1) {};
\node[cell] (b2) at (1,-1) {};
\node[cell] (b3) at (2,-1) {};
\node[cell] (b4) at (3,-1) {};
\node[cell] (b5) at (4,-1) {};
\node[cell] (b6) at (5,-1) {};

\draw[wall] (0.5,0.5) -- (0.5,-0.5);
\draw[wall] (0.5,-0.5) -- (1.5,-0.5);
\draw[wall] (2.5,0.5) -- (2.5,-0.5);
\draw[wall] (2.5,-0.5) -- (3.5,-0.5);
\draw[wall] (4.5,0.5) -- (4.5,-0.5);
\draw[wall] (4.5,-0.5) -- (5.5,-0.5);
\end{tikzpicture}}
\caption{A Young building $\mathcal{B}$ (denoted by Z3) on the left and a Young building $\mathcal{B}^3=(\mathrm{Z}3)^3$ with periodic walls on the right.}\label{tab:youngZ3}
\end{figure}

Banderier and Wallner \cite{banderier2021young} counted Young tableaux with periodic walls by the density method and bijection with binary trees. Let $\lambda=(2,2)$.
They considered a classification of $2\times 2$ periodic shapes (a total of $2^6=64$ different models).
In other words, there are $64$ distinct sets $\mathcal{W}$.
Two of them are shown in \cref{tab:youngZ2,tab:youngZ3}.
From the perspective of counting $|\mathrm{YT}(\mathcal{B}^n)|$, there are actually $32$ models in total by symmetry.
Based on their counting formulas, Banderier and Wallner classified these $32$ counting models into the following four classes:
(P) simple products, (A) algebraic, (H) hypergeometric, (Z) zig-zag.
They also obtained exact counting formulas for 30 of these 32 counting models; see \cite[Table 1]{banderier2021young}.
The formulas for the remaining counting models in \cref{tab:youngZ2,tab:youngZ3} are unknown.
They are denoted by Z2 and Z3, respectively.

\begin{figure}[htbp]
\centering
\subfloat[]{
\begin{tikzpicture}[
    cell/.style={draw, minimum size=1cm, anchor=center},
    wall/.style={line width=2pt, red}]
\node[cell] (a1) at (0,0) {};
\node[cell] (a2) at (1,0) {};

\node[cell] (b1) at (0,-1) {};
\node[cell] (b2) at (1,-1) {};

\draw[wall] (-0.5,0.5) -- (-0.5,-0.5);
\draw[wall] (-0.5,-0.5) -- (0.5,-0.5);
\draw[wall] (1.5,0.5) -- (1.5,-0.5);
\end{tikzpicture}}
\hspace{15pt}
\subfloat[]{
\begin{tikzpicture}[
    cell/.style={draw, minimum size=1cm, anchor=center},
    wall/.style={line width=2pt, red}]
\node[cell] (a1) at (0,0) {};
\node[cell] (a2) at (1,0) {};
\node[cell] (a3) at (2,0) {};
\node[cell] (a4) at (3,0) {};
\node[cell] (a5) at (4,0) {};
\node[cell] (a6) at (5,0) {};

\node[cell] (b1) at (0,-1) {};
\node[cell] (b2) at (1,-1) {};
\node[cell] (b3) at (2,-1) {};
\node[cell] (b4) at (3,-1) {};
\node[cell] (b5) at (4,-1) {};
\node[cell] (b6) at (5,-1) {};

\draw[wall] (-0.5,0.5) -- (-0.5,-0.5);
\draw[wall] (-0.5,-0.5) -- (0.5,-0.5);
\draw[wall] (1.5,0.5) -- (1.5,-0.5);
\draw[wall] (1.5,-0.5) -- (2.5,-0.5);
\draw[wall] (3.5,0.5) -- (3.5,-0.5);
\draw[wall] (3.5,-0.5) -- (4.5,-0.5);
\draw[wall] (5.5,0.5) -- (5.5,-0.5);
\end{tikzpicture}}
\caption{A Young building $\mathcal{B}$ (denoted by Z2) on the left and a Young building $\mathcal{B}^3=(\mathrm{Z}2)^3$ with periodic walls on the right.}\label{tab:youngZ2}
\end{figure}

Let $a_n$ and $b_n$ count the Young tableaux with periodic walls of types Z2 and Z3 with $n$ blocks.
That is
$$a_n=|\mathrm{YT}((\mathrm{Z2})^n)|\qquad\text{and}\qquad b_n=|\mathrm{YT}((\mathrm{Z3})^n)|.$$
Each such tableau has $4n$ cells, and the empty tableau has counts $a_0 = b_0 = 1$.
We first specify the normalization, since the number of blocks and the number of cells give different series.
Define
\begin{equation}\label{eq:defAB}
A(z)=\sum_{n \geq 0} \frac{a_{n}}{(4 n) !} z^{n},\qquad
B(z)=\sum_{n \geq 0} \frac{b_{n}}{(4 n) !} z^{n}.
\end{equation}
The exponential generating functions in which the variable marks cells are therefore $A(t^{4})$ and $B(t^{4})$. The usual exponential generating series of the sequences indexed by the number of blocks are
\begin{equation}\label{eq:defF2F3}
F_{2}(t)=\sum_{n \geq 0} \frac{a_{n}}{n !} t^{n},\qquad
F_{3}(t)=\sum_{n \geq 0} \frac{b_{n}}{n !} t^{n}.
\end{equation}

In particular, Banderier and Wallner proposed the following conjecture.

\begin{conjecture}{(Banderier and Wallner~\cite{banderier2021young})}\label{MainConjecture-Z2Z3}
The exponential generating functions $F_2(t)$ and $F_3(t)$ are D-algebraic but not D-finite.
\end{conjecture}

The main motivation of this paper is to resolve \cref{MainConjecture-Z2Z3}.
We give counting formulas for both families and prove that their exponential generating functions are in fact differentially transcendental. These series require a separate analysis: changing a factorial normalization does not in general preserve differential algebraicity.

A function $F(z)$ is \emph{differentially algebraic} (or \emph{D-algebraic}) over $\mathbb{C}(z)$ if it satisfies
$$P\big(z, F(z), F'(z), \dots, F^{(k)}(z)\big)=0$$
for some integer $k\geq 0$ and some nonzero polynomial $P\in\mathbb{C}(z)[Y_0,\dots,Y_k]$.
It is \emph{differentially transcendental} if no such polynomial exists.
It is \emph{D-finite} if it satisfies a nonzero linear differential equation with coefficients in $\mathbb{C}(z)$.
These definitions also apply to formal series.
A sequence $(u_n)_{n\ge 0}$ is \emph{P-recursive} if, for some $d\ge 0$, there are polynomials $p_0,\dots,p_d\in\mathbb{C}[n]$, not all zero, such that
\[
\sum_{j=0}^{d} p_{j}(n)\, u_{n+j}=0
\]
for every sufficiently large integer $n$. Equivalently, its ordinary generating series is D-finite; see~\cite[Theorem~1.5]{Stanley80}.
For further details on the definitions, we refer the reader to \cite[Chapters 5 and 6]{Stanley99}.

\begin{theorem}\label{thm:ABnotDfinite}
The generating functions $A(z)$ and $B(z)$ are defined in \eqref{eq:defAB}.
Then neither $A(z)$ nor $B(z)$ is D-finite over $\mathbb{C}(z)$. Equivalently, the sequences $(a_n)_{n\ge0}$ and $(b_n)_{n\ge0}$ are not P-recursive. Moreover,
\begin{equation*}
\lim_{n\to\infty}\left(\frac{(2n)!}{(4n)!}a_n\right)^{1/n}
=
\lim_{n\to\infty}\left(\frac{(2n)!}{(4n)!}b_n\right)^{1/n}
=\frac{\pi^2}{32},
\end{equation*}
and the associated auxiliary series has radius of convergence $32/\pi^2$, which is not an algebraic number.
\end{theorem}

\begin{theorem}\label{thm:block-egfs}
The series $F_2(t)$ and $F_3(t)$ in \eqref{eq:defF2F3} have radius of convergence zero. As formal power series, each is differentially transcendental over $\mathbb{C}(t)$. Consequently, neither series is D-algebraic, and neither is D-finite.
More precisely,
\[
\lim _{n \to \infty} \frac{a_{n+1}}{n(n+1) a_{n}}=\lim _{n \to \infty} \frac{b_{n+1}}{n(n+1) b_{n}}=2 \pi^{2}
\]
and
\[
\lim _{n\to\infty}\left(\frac{a_{n}}{(n!)^{2}}\right)^{1/n}=\lim _{n \to \infty}\left(\frac{b_{n}}{(n !)^{2}}\right)^{1/n}=2\pi^{2}.
\]
\end{theorem}

We also obtained two recurrences for $a_n$ and $b_n$; see \cref{thm:counts}.
This allows us to compute certain concrete values quickly.
Computing all counts with $0\leq n\leq N$ requires $O(N^2)$ arithmetic operations and $O(N)$ entries of working memory.
As a byproduct, non-recursive formulas for $a_n$ and $b_n$ are given in \cref{prop:subset-sums}.
However, they are not convenient for computation.

The rest of this paper is organized as follows.
In \cref{sec:counts}, we associate Z2 and Z3 with their corresponding posets, and then obtain recurrences for $a_n$ and $b_n$.
\cref{sec:growth} is mainly devoted to proving \cref{thm:ABnotDfinite}.
\cref{sec:block-egfs} mainly concerns the proof of \cref{thm:block-egfs}.
\cref{sec:subset-sums} is a byproduct; it mainly contains the non-recursive formulas for $a_n$ and $b_n$.

\section{Two posets and their counting formulas}\label{sec:counts}

This section gives the details of the two families of posets, denoted by
$\Ztwo$ and $\Zthree$, and the counts $a_n$ and $b_n$ of their linear
extensions with $n$ blocks. These objects were introduced in the
Introduction. We describe the order relations in each family. We then
express $a_n$ and $b_n$ as integrals over order polytopes. These integral
formulas are the basis for the later analysis. We also record a polynomial
recursion, which gives single sums with integer coefficients. This
recursion is useful for computation, but the later arguments do not rely
on it.

\subsection{Order relations and order polytopes}

We describe $\Ztwo$ and $\Zthree$ by their order relations.
In this way, the objects do not depend on a diagram convention.

Write $x_1,\ldots,x_{2n}$ for the bottom-row cells and
$y_1,\ldots,y_{2n}$ for the top-row cells, read from left to right.
In both families the bottom row is a chain:
\[
 x_1<x_2<\cdots<x_{2n}.
\]
For $\Ztwo$, the additional generating relations are
\begin{equation}
 y_{2j-1}<y_{2j},\qquad x_{2j}<y_{2j} \quad(1\leq j\leq n).
 \label{eq:Z2-rel}
\end{equation}
For $\Zthree$, they are
\begin{equation}
 x_{2j-1}<y_{2j-1}\quad(1\leq j\leq n),
 \qquad
 y_{2j}<y_{2j+1}\quad(1\leq j<n).
 \label{eq:Z3-rel}
\end{equation}
The orders are the transitive closures of the stated relations. No
other adjacent-cell relations are imposed. In particular, in
\eqref{eq:Z3-rel}, $y_{2n}$ is isolated. Also, $y_1$ has only its
attachment to $x_1$.

A legal labeling is a bijection from the cells to
$\{1,\ldots,4n\}$ that preserves these orders. Thus it is a linear
extension of the corresponding poset.

\begin{figure}[ht]
\centering
\begin{tikzpicture}[x=1.12cm,y=1cm,>=stealth]
\foreach \j in {1,2,3,4} {
  \node (a\j) at (\j,0) {$x_{\j}$};
  \node (b\j) at (\j,1) {$y_{\j}$};
}
\draw[->] (a1)--(a2); \draw[->] (a2)--(a3); \draw[->] (a3)--(a4);
\draw[->] (b1)--(b2); \draw[->] (b3)--(b4);
\draw[->] (a2)--(b2); \draw[->] (a4)--(b4);
\node at (2.5,-.65) {$\Ztwo$};
\end{tikzpicture}
\qquad
\begin{tikzpicture}[x=1.12cm,y=1cm,>=stealth]
\foreach \j in {1,2,3,4} {
  \node (a\j) at (\j,0) {$x_{\j}$};
  \node (b\j) at (\j,1) {$y_{\j}$};
}
\draw[->] (a1)--(a2); \draw[->] (a2)--(a3); \draw[->] (a3)--(a4);
\draw[->] (a1)--(b1); \draw[->] (a3)--(b3); \draw[->] (b2)--(b3);
\node at (2.5,-.65) {$\Zthree$};
\end{tikzpicture}
\caption{The two posets with two blocks. Only generating order
relations are drawn.}
\label{fig:posets}
\end{figure}
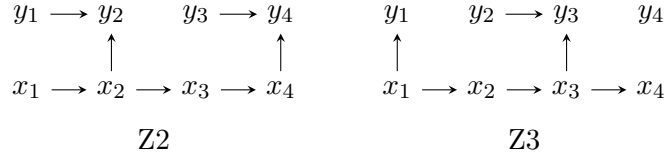

\cref{fig:posets} shows the two orders for $n=2$.
Every arrow points from a smaller element to a larger one.
The isolated element in the second diagram is part of the
counting problem and contributes to the factorial $(4n)!$.

For a finite poset $P$ with $N$ elements, its order polytope is
\[
 \mathcal O(P)=\{(u_v)_{v\in P}\in[0,1]^N:
                 u_v\leq u_w\ \hbox{whenever }v\leq_Pw\}.
\]
The volume formula
\begin{equation}
 \#\{\hbox{linear extensions of }P\}=N!\vol\mathcal O(P)
 \label{eq:order-volume}
\end{equation}
is due to Stanley~\cite[Corollary~4.2]{Stanley86}. For completeness,
the hyperplanes on which two coordinates coincide have measure zero.
Off these hyperplanes, each ordering of all $N$ coordinates describes
an open simplex of volume $1/N!$. Such a simplex belongs to
$\mathcal O(P)$ exactly when the ordering is a linear extension. These
simplices have disjoint interiors. They cover the order polytope up to
measure zero. This proves \eqref{eq:order-volume}.

Set
\[
 q(x)=\frac{1-x^2}{2},\qquad 0\leq x\leq1.
\]
We use the same symbols for cells and their coordinates in the order
polytope. Once the bottom-row coordinates are fixed, the coordinates
in different top-row components can be integrated independently.

\begin{proposition}
For $n\geq1$, we have 
\begin{align}
 \frac{a_n}{(4n)!}
 &=\int_{0<x_1<\cdots<x_{2n}<1}
       \prod_{j=1}^n q(x_{2j})\dd x_1\cdots\dd x_{2n},
 \label{eq:a-density}\\
 \frac{b_n}{(4n)!}
 &=\int_{0<x_1<\cdots<x_{2n}<1}
       (1-x_1)\prod_{j=2}^n q(x_{2j-1})
       \dd x_1\cdots\dd x_{2n}.
 \label{eq:b-density}
\end{align}
\end{proposition}

\begin{proof}
For \eqref{eq:Z2-rel}, fix the bottom-row coordinates. For each $j$, the
coordinates $y_{2j-1}$ and $y_{2j}$ only interact with $x_{2j}$. The
relations are $y_{2j-1}<y_{2j}$ and $x_{2j}<y_{2j}$. Thus, for fixed
$x_{2j}$, we have $x_{2j}<y_{2j}<1$ and $0<y_{2j-1}<y_{2j}$. Their
integral is
\[
 \int_{x_{2j}}^1 \int_0^{y_{2j}} dy_{2j-1} dy_{2j}
 = \int_{x_{2j}}^1 y_{2j} dy_{2j}
 = q(x_{2j}).
\]
The pairs for different $j$ are independent. Multiplying over all $j$ and
then integrating over the bottom chain gives \eqref{eq:a-density} by
\eqref{eq:order-volume}.

For \eqref{eq:Z3-rel}, again fix the bottom-row coordinates. The
coordinate $y_1$ only interacts with $x_1$. The relation is $x_1<y_1$.
Thus its integral is
\[
 \int_{x_1}^1 dy_1 = 1-x_1.
\]
For $2\leq j\leq n$, the pair $(y_{2j-2},y_{2j-1})$ only interacts with
$x_{2j-1}$. The relations are $y_{2j-2}<y_{2j-1}$ and
$x_{2j-1}<y_{2j-1}$. For fixed $x_{2j-1}$, the integral is
\[
 \int_{x_{2j-1}}^1 \int_0^{y_{2j-1}} dy_{2j-2} dy_{2j-1}
 = \int_{x_{2j-1}}^1 y_{2j-1} dy_{2j-1}
 = q(x_{2j-1}).
\]
The isolated coordinate $y_{2n}$ contributes a factor $1$. These are all
the top-row components. Multiplying them and integrating over the bottom
chain gives \eqref{eq:b-density}.
\end{proof}

\subsection{A polynomial recursion}

To write the density formulas in \eqref{eq:a-density} and
\eqref{eq:b-density} in iterative form, we introduce an integral
operator. Define
\begin{equation}
 (Tf)(x)=\int_0^x(x-u)q(u)f(u)\dd u.
 \label{eq:T}
\end{equation}
For $k\geq0$, direct integration gives
\begin{equation}
 T(x^k)=\frac{x^{k+2}}{2(k+1)(k+2)}
        -\frac{x^{k+4}}{2(k+3)(k+4)}.
 \label{eq:T-monomial}
\end{equation}
Thus $T$ maps $\Q[x]$ into itself.

In the following proposition, subscripts denote partial derivatives with
respect to $x$. Thus $S_x=\partial S/\partial x$ and
$S_{xx}=\partial^2 S/\partial x^2$. The same convention applies to $V$.

\begin{proposition}
The counts can be evaluated by
\begin{equation}
 \frac{a_n}{(4n)!}=\left.\frac{\partial}{\partial x}(T^n x)\right|_{x=1},
 \qquad
 \frac{b_n}{(4n)!}
 =\left.T^{n-1}\left(\frac{x^2}{2}-\frac{x^3}{6}\right)\right|_{x=1}
 \quad(n\geq1).
 \label{eq:operator-counts}
\end{equation}
Let $S,V$ be the solutions of
\begin{align}
 S_{xx}(x,z)&=zq(x)S(x,z),&S(0,z)&=0,&S_x(0,z)&=1,
 \label{eq:S-ivp}\\
 V_{xx}(x,z)&=zq(x)V(x,z)+z(1-x),&V(0,z)&=0,&V_x(0,z)&=0.
 \label{eq:V-ivp}
\end{align}
Then
\begin{equation}
 A(z)=S_x(1,z),\qquad B(z)=1+V(1,z).
 \label{eq:boundary-values}
\end{equation}
The functions in \eqref{eq:boundary-values} are entire in $z$.
\end{proposition}

\begin{proof}
We first record two properties of $T$. From \eqref{eq:T}, we have 
\[
 (Tf)(0)=0,\qquad (Tf)'(0)=0,\qquad (Tf)''(x)=q(x)f(x).
\]
We use these properties below.

The polynomial $T^n x$ is the integral over
$0<u_1<\cdots<u_{2n+1}<x$ with weight
$\prod_{j=1}^nq(u_{2j})$. This statement follows by induction.
For $n=0$, the integral is $\int_0^x\dd u_1=x$. To pass from $n-1$ to
$n$, applying $T$ adds the coordinate $u_{2n}$ with weight
$q(u_{2n})$. The factor $x-u_{2n}$ can be rewritten as the integral over
the last coordinate $u_{2n+1}$. Thus
\[
 T^n x
 =\int_{0<u_1<\cdots<u_{2n+1}<x}
       \prod_{j=1}^nq(u_{2j})\dd u_1\cdots\dd u_{2n+1}.
\]
Differentiating this integral with respect to the upper limit $x$
removes the last integration. By the Leibniz rule, the derivative is
\[
 \frac{\partial}{\partial x}(T^n x)
 =\int_{0<u_1<\cdots<u_{2n}<x}
       \prod_{j=1}^nq(u_{2j})\dd u_1\cdots\dd u_{2n}.
\]
At $x=1$, the result is \eqref{eq:a-density}.

Similarly,
\[
 \frac{x^2}{2}-\frac{x^3}{6}
 =\int_{0<u_1<u_2<x}(1-u_1)\dd u_1\dd u_2.
\]
Each further application of $T$ appends two coordinates and puts a
factor $q$ at the first appended coordinate. The weight obtained in
this way is exactly the integrand in \eqref{eq:b-density}. Hence
\eqref{eq:operator-counts} follows.

The relevant series are
\begin{align*}
 S(x,z)=\sum_{n\geq0}z^nT^n x\qquad\text{and}\qquad
 V(x,z)=\sum_{n\geq1}z^nT^{n-1}\left(\frac{x^2}{2}-\frac{x^3}{6}\right).
\end{align*}
They satisfy $S=x+zTS$, and
$V=z(x^2/2-x^3/6)+zTV$. Differentiating twice in $x$ gives
\eqref{eq:S-ivp}--\eqref{eq:V-ivp}, including the initial values.
Equation~\eqref{eq:operator-counts} gives
\eqref{eq:boundary-values}.

For convergence, we use the multiple-integral form of $T^n$. Since
$0\leq q(u)\leq1/2$ on $[0,1]$, we obtain the following estimate. For
$f\in C([0,1])$ and $n\ge 0$, set
\[
 \|f\|_{\infty}=\sup_{x\in[0,1]}|f(x)|.
\]
Then
\[
\begin{aligned}
 |(T^n f)(x)|
 &=\left|\int_{0<u_1<\cdots<u_{2n}<x}
   \prod_{j=0}^{n-1} q(u_{2j+1})\, f(u_1)
   \,du_1\cdots du_{2n}\right|\\
 &\le \|f\|_{\infty} \cdot \frac{1}{2^n}
   \int_{0<u_1<\cdots<u_{2n}<x}
   du_1\cdots du_{2n}\\
 &=\|f\|_{\infty}\frac{x^{2n}}{2^n\,(2n)!},
 \qquad 0\le x\le 1.
\end{aligned}
\]
Hence the displayed series converge locally uniformly in $z$. The
convergence is uniform for $x\in[0,1]$. Their first $x$-derivatives
have analogous bounds obtained by differentiating \eqref{eq:T}. Their
second derivatives are controlled by the integral equations. This
justifies the differentiations. Alternatively, the density integrals
give
\[
 0<\frac{a_n}{(4n)!}\leq\frac{1}{2^n(2n)!},\qquad
 0<\frac{b_n}{(4n)!}\leq\frac{1}{2^{n-1}(2n)!}\quad(n\geq1).
\]
These inequalities prove directly that $A$ and $B$ are entire.
\end{proof}

\subsection{Single sums with integer coefficients}

The operator formula already gives a short exact calculation.
The following form avoids polynomial integration during the
calculation. It uses integer coefficients throughout the recurrence.
For a positive integer parameter $s$, define
\begin{equation} \label{eq:integer-triangle}
 c^{(s)}_{0,0}=1,\qquad
 c^{(s)}_{n+1,k}=c^{(s)}_{n,k}+(s+2n+2k-1)(s+2n+2k)c^{(s)}_{n,k-1}.
\end{equation}
Here $c^{(s)}_{n,k}=0$ when $k<0$ or $k>n$. The recurrence is
used for $n\geq0$ and $0\leq k\leq n+1$.
For each positive integer $s$, every nonzero entry is a positive
integer. In this paper we only need the cases $s=1,2,3$. The initial
polynomials in \eqref{eq:operator-counts} are linear combinations of
$x^s/s!$ for $s=1,2,3$.

\begin{theorem}\label{thm:counts}
For $n\geq0$, we have 
\begin{equation}\label{eq:a-single}
 a_n=\frac{(4n)!}{2^n} \sum_{k=0}^n\frac{(-1)^k c^{(1)}_{n,k}}{(2n+2k)!}.
\end{equation}
For $n\geq1$, we have
\begin{equation} \label{eq:b-single}
 b_n=\frac{(4n)!}{2^{n-1}}\sum_{k=0}^{n-1}(-1)^k\left(\frac{c^{(2)}_{n-1,k}}{(2n+2k)!}-\frac{c^{(3)}_{n-1,k}}{(2n+2k+1)!}\right).
\end{equation}
The coefficients are determined by \eqref{eq:integer-triangle}.
Computing all counts with $0\leq n\leq N$ requires $O(N^2)$ arithmetic operations and $O(N)$ entries of working memory.
\end{theorem}

\begin{proof}
The factorial basis makes the calculation explicit. For $m\geq0$,
\eqref{eq:T-monomial} becomes
\[
 2T\left(\frac{x^m}{m!}\right)
 =\frac{x^{m+2}}{(m+2)!}
  -(m+1)(m+2)\frac{x^{m+4}}{(m+4)!}.
\]
We claim that
\begin{equation}
 2^nT^n\left(\frac{x^s}{s!}\right)
 =\sum_{k=0}^n(-1)^kc^{(s)}_{n,k}
       \frac{x^{s+2n+2k}}{(s+2n+2k)!}.
 \label{eq:factorial-basis}
\end{equation}
At $n=0$, the claim is the trivial identity $x^s/s!=x^s/s!$.
Assume it holds at $n$. Applying $2T$, the first term in the
preceding monomial formula contributes $c^{(s)}_{n,k}$ to the
coefficient with index $k$ at level $n+1$. The second term with
old index $k-1$ has degree $m=s+2n+2k-2$ before applying $T$.
Its extra minus sign changes $(-1)^{k-1}$ to $(-1)^k$. Its
integer multiplier is
\[
 (m+1)(m+2)=(s+2n+2k-1)(s+2n+2k).
\]
There are no other contributions to this degree. Thus its total
coefficient is exactly \eqref{eq:integer-triangle}. This proves
\eqref{eq:factorial-basis} by induction.

For $a_n$, use $s=1$ in \eqref{eq:factorial-basis}, differentiate
in $x$, and evaluate at $x=1$. The differentiation changes
$x^{2n+2k+1}/(2n+2k+1)!$ into
$x^{2n+2k}/(2n+2k)!$. The first identity in
\eqref{eq:operator-counts} gives \eqref{eq:a-single}.

For $b_n$, write the initial polynomial as
$x^2/2!-x^3/3!$. Use \eqref{eq:factorial-basis} with $n-1$
in place of $n$, first with $s=2$ and then with $s=3$.
Evaluation at $x=1$ and the second identity in
\eqref{eq:operator-counts} give \eqref{eq:b-single}.

There are $n+1$ entries in row $n$ of each triangle. Each
entry requires a bounded number of integer additions and
multiplications. Hence the rows through $N$ require $O(N^2)$
arithmetic operations. Consecutive factorials up to $(4N+1)!$
can be computed in $O(N)$ further operations. The displayed
sums use $O(n)$ rational operations at index $n$. Therefore
they use $O(N^2)$ operations in total. Each new row needs only the
preceding row. Retaining the current rows and factorials uses
$O(N)$ entries. No bound on the number of bits in an entry is
asserted.
\end{proof}

For example, the first three rows for $s=1$ are
\[
 (c^{(1)}_{0,k})=(1),\qquad
 (c^{(1)}_{1,k})=(1,6),\qquad
 (c^{(1)}_{2,k})=(1,26,252).
\]
Thus \eqref{eq:a-single} gives
\[
 a_1=\frac{4!}{2}\left(\frac1{2!}-\frac6{4!}\right)=3,
 \qquad
 a_2=\frac{8!}{4}\left(\frac1{4!}-\frac{26}{6!}
                                      +\frac{252}{8!}\right)=119.
\]
For $b_1$, the two initial entries $c^{(2)}_{0,0}=c^{(3)}_{0,0}=1$
give $b_1=4!(1/2!-1/3!)=8$.

Expanding the integrands directly also gives the subset sums
in \cref{sec:subset-sums}. Those formulas contain $2^n$
terms; \eqref{eq:a-single} and \eqref{eq:b-single} collect the
calculation into rows of the integer recurrence.

The formulas yield the following values. They also agree with
direct enumeration of linear extensions for all entries displayed.
\begin{center}
\begin{tabular}{@{}rrr@{}}
\toprule
$n$ & $a_n$ & $b_n$\\
\midrule
0 & 1 & 1\\
1 & 3 & 8\\
2 & 119 & 416\\
3 & 13\,761 & 56\,136\\
4 & 3\,178\,785 & 14\,433\,600\\
5 & 1\,226\,341\,035 & 6\,042\,488\,040\\
6 & 711\,310\,157\,271 & 3\,743\,684\,613\,216\\
\bottomrule
\end{tabular}
\end{center}

\section{The elementary proof of \cref{thm:ABnotDfinite}}\label{sec:growth}

In this section we prove \cref{thm:ABnotDfinite}. The theorem
states the exponential growth of the normalized counts and the
non-$P$-recursiveness of $(a_n)$ and $(b_n)$. Our proof uses the
integral formulas from \cref{sec:counts} and a standard fact
about $D$-finite series. We begin with the normalization.

Let
\[
 L=\int_0^1\sqrt{q(x)}\dd x=\frac{\pi}{4\sqrt2},\qquad
 d_n=\frac{(2n)!}{(4n)!}a_n,\qquad
 e_n=\frac{(2n)!}{(4n)!}b_n.
\]
We first prove the limits
\[
 \lim_{n\to\infty}d_n^{1/n}
 =\lim_{n\to\infty}e_n^{1/n}
 =L^2=\frac{\pi^2}{32}.
\]

Since $q$ decreases on $[0,1]$, we have
$q(x_{2j})\leq\sqrt{q(x_{2j-1})q(x_{2j})}$. The product
$\prod_{k=1}^{2n}\sqrt{q(x_k)}$ is symmetric. Its integral over the
ordered simplex is therefore $1/(2n)!$ times its integral over the
cube. By \eqref{eq:a-density},
\begin{equation}
 d_n\leq L^{2n}.
 \label{eq:d-upper}
\end{equation}

For a lower bound, fix \(c\in(0,1)\) and restrict the integral to
\(x_{2n}<c\). This avoids the zero of \(q\) at \(1\). On that simplex,
\[
 \prod_{j=1}^n\frac{q(x_{2j})}{q(x_{2j-1})}
 =\frac{q(x_{2n})}{q(x_1)}
  \prod_{j=1}^{n-1}\frac{q(x_{2j})}{q(x_{2j+1})}
 \geq\frac{q(c)}{q(0)}.
\]
It follows that
\begin{equation}
 d_n\geq\sqrt{\frac{q(c)}{q(0)}}
             \left(\int_0^c\sqrt{q(x)}\dd x\right)^{2n}.
 \label{eq:d-lower}
\end{equation}
Taking the \(n\)-th roots in \eqref{eq:d-upper} and \eqref{eq:d-lower} gives
\[
 \limsup_{n\to\infty} d_n^{1/n}\le L^2
 \qquad\text{and}\qquad
 \liminf_{n\to\infty} d_n^{1/n}\ge \sup_{0<c<1}\left(\int_0^c\sqrt{q(x)}\,dx\right)^2=L^2.
\]
Hence \(d_n^{1/n}\to L^2\).

For the second family, we compare the integrands. Since
\(1-x_1\ge q(x_1)\ge q(x_2)\) and
\(q(x_{2j-1})\ge q(x_{2j})\) for \(j\ge2\), the integrand in
\eqref{eq:b-density} is at least the integrand in
\eqref{eq:a-density}. Hence \(e_n\ge d_n\). Consequently,
\[
 \liminf_{n\to\infty} e_n^{1/n}
 \ge \lim_{n\to\infty} d_n^{1/n}=L^2.
\]

For the reverse inequality, use \(1-x_1\le1\) and, for \(j\ge2\),
\[
 q(x_{2j-1})\le\sqrt{q(x_{2j-2})q(x_{2j-1})}.
\]
Therefore
\[
 (1-x_1)\prod_{j=2}^n q(x_{2j-1})
 \le\prod_{k=2}^{2n-1}\sqrt{q(x_k)}.
\]
For \(n\ge2\), integrate \(x_1\) and \(x_{2n}\) first. The resulting
lengths \(x_2\) and \(1-x_{2n-1}\) are at most \(1\). Thus
\[
 \frac{b_n}{(4n)!}\le\frac{L^{2n-2}}{(2n-2)!}.
\]
For \(n=1\), the same inequality holds because
\(b_1/4!=1/3\le1\). Hence
\[
 d_n\le e_n\le(2n)(2n-1)L^{2n-2}.
\]
Taking \(n\)th roots gives
\[
 \limsup_{n\to\infty} e_n^{1/n}\le L^2.
\]
Together with the lower bound, this yields
\(e_n^{1/n}\to L^2\). In particular,
\[
 a_n=\frac{(4n)!}{(2n)!}
          \left(\frac{\pi^2}{32}\right)^n\ee^{o(n)},\qquad
 b_n=\frac{(4n)!}{(2n)!}
          \left(\frac{\pi^2}{32}\right)^n\ee^{o(n)}.
\]

Next we show that the generating functions are not $D$-finite. Recall that a formal power series is
$D$-finite if it satisfies a linear differential equation with polynomial coefficients. We use the following
standard fact. See, for instance, \cite[Chapter~6]{Stanley99} or \cite{Lipshitz89}.

\begin{lemma}
\label{lem:radius}
A $D$-finite series with rational coefficients and finite positive
radius of convergence has an algebraic radius of convergence.
\end{lemma}

\begin{proof}
First, an equation with coefficients in $\C[z]$ can be replaced by
one with coefficients in $\Q[z]$. To see this, fix the order and
degree bounds of the given equation, and regard its finitely many
coefficients as unknowns. Equating coefficients of powers of $z$
gives a homogeneous linear system over $\Q$. The row space of this
possibly infinite system is spanned by finitely many rows. Its
kernel has the same dimension over $\Q$ as after scalar extension
to $\C$. Thus a nonzero complex solution implies a nonzero rational
solution.

We now use the equation with rational coefficients. Let
\[
 p_m(z)F^{(m)}(z)+\cdots+p_0(z)F(z)=0,
 \qquad p_j\in\Q[z],\quad p_m\ne0,
\]
be such an equation. A germ of a solution continues analytically
through an ordinary point. Indeed, after division by $p_m$, this
follows from the existence theorem for an analytic first-order
linear system applied to $(F,\ldots,F^{(m-1)})$. Hence every
singular point of the germ on its finite circle of convergence is
a zero of $p_m$. At least one such singular point exists; otherwise
finitely many neighborhoods covering that circle would extend $F$
to a larger disk. Its radius is therefore $|\alpha|$ for a zero
$\alpha$ of $p_m$. Both $\alpha$ and its complex conjugate are
algebraic, and $|\alpha|^2=\alpha\overline\alpha$. Thus $|\alpha|$
is algebraic.
\end{proof}

The ordinary series $\sum d_nz^n$ and $\sum e_nz^n$ have radius
$L^{-2}=32/\pi^2$, because $d_n^{1/n}\to L^2$ and
$e_n^{1/n}\to L^2$. This number is transcendental by the
transcendence of $\pi$~\cite{Lindemann82}. By \cref{lem:radius},
neither series is $D$-finite. Thus neither $(d_n)$ nor $(e_n)$ is
$P$-recursive.

For each fixed $j\geq0$, the ratio
\[
 \frac{(2n+2j)!/(4n+4j)!}{(2n)!/(4n)!}
\]
is a rational function of $n$. If a polynomial-coefficient recurrence
existed for $(a_n)$, substitute $a_n=d_n(4n)!/(2n)!$ and clear
denominators. The result would be a nonzero polynomial-coefficient
recurrence for $(d_n)$, a contradiction. The same argument applies to
$(b_n)$ and $(e_n)$. Using the factorial factor $(4n)!$ instead, we
also find that neither $a_n/(4n)!$ nor $b_n/(4n)!$ is $P$-recursive.
Hence $A$ and $B$ are not $D$-finite. This proves \cref{thm:ABnotDfinite}.

\cref{sec:growth} proved that \(A\) and \(B\) are not
\(D\)-finite. The next section disproves 
the Banderier--Wallner conjecture. That is, we show that 
the block-counted series \(F_2\) and \(F_3\) are differentially
transcendental. The method is different: a coefficient criterion
based on ratio limits, combined with a probabilistic representation
of the normalized counts.
\section{Disproving the Banderier–Wallner conjecture}
\label{sec:block-egfs}

We prove \cref{thm:block-egfs} for the block-counted
exponential generating functions $F_2(t),F_3(t)$ in
\eqref{eq:defF2F3}. Both are differentially transcendental
as formal series, and their coefficient ratios tend to
$2\pi^2$ after division by $n$. Their radii of convergence are
zero, while their ordinary Borel transforms, obtained by
dividing the coefficient of $t^n$ by one further factor $n!$,
have radius $1/(2\pi^2)$. The differential-algebraic classification
of $F_2,F_3$ is understood coefficient by coefficient in $\C[[t]]$.

Radius zero alone does not imply differential transcendence. For
example, the Euler series $E(t)=\sum_{n\geq0}n!t^n$ has radius
zero and satisfies
\[
 t^2E'(t)+(t-1)E(t)+1=0.
\]
Differentiation gives the homogeneous linear equation
\[
 t^2E''(t)+(3t-1)E'(t)+E(t)=0,
\]
so $E$ is $D$-finite. The obstruction used here is the
transcendental value of the coefficient ratio, together with
rationality and positivity of the coefficients.

In the Euler example, \(c_n=n!\) and \(c_{n+1}/(nc_n)=(n+1)/n\to1\).
The limit \(1\) is algebraic, so \cref{thm:ratio-algebraic} below does not exclude
differential algebraicity. The obstruction in our application is
that the limit \(2\pi^2\) is transcendental.

\subsection{A coefficient criterion for differential transcendence}

We first prove the factorial convolution estimates needed for the
criterion.

\begin{lemma}
\label{lem:factorial-convolution}
Let $(u_n)_{n\geq0}$ be a sequence of positive real numbers such that
\[
 \lim_{n\to\infty}\frac{u_{n+1}}{nu_n}=K,\qquad 0<K<\infty.
\]
For every fixed nonnegative integer $M$, we have 
\[
 \sum_{j=M}^{n-M}u_j u_{n-j}=O(u_n n^{-M}).
\]
For every fixed integer $q\geq2$, the corresponding restricted $q$-fold convolution also satisfies
\[
\sum_{\substack{n_1+\cdots+n_q=n\\n_1,\ldots,n_q\geq M}}u_{n_1}\cdots u_{n_q}=O(u_n n^{-M}).
\]
The constants in these estimates may depend on $M$ and $q$.
\end{lemma}

\begin{proof}
Choose $\delta\in(0,1)$ sufficiently small that
\[
 \gamma=\frac{1+\delta}{1-\delta}<2.
\]
Since $n/(n+1)\to1$, there is an integer $n_0$ such that
\[
 K(1-\delta)(j+1)u_j
 \leq u_{j+1}
 \leq K(1+\delta)(j+1)u_j
 \qquad(j\geq n_0).
\]

After enlarging the constant to cover the finitely many indices
$0\leq j<n_0$, the upper inequality gives
\[
 u_j\leq C K^j(1+\delta)^j j!
 \qquad(j\geq0).
\]
If $0\leq j\leq n/2$ and $n\geq2n_0$, telescoping the lower
inequality from $n-j$ to $n$ gives
\[
 \frac{u_{n-j}}{u_n}
 \leq\frac{(n-j)!}{n!}\,[K(1-\delta)]^{-j}.
\]
Combining these bounds yields
\[
 \frac{u_j u_{n-j}}{u_n}
 \leq C\frac{\gamma^j}{\binom nj}
 \qquad(0\leq j\leq n/2).
\]

We estimate the sum on the right beginning at any fixed $M$.
For $M\leq j\leq n/4$, the ratio of consecutive terms is
\[
 \frac{\gamma^{j+1}/\binom n{j+1}}
      {\gamma^j/\binom nj}
 =\gamma\frac{j+1}{n-j}.
\]
For all sufficiently large $n$, this ratio is bounded by a fixed
number less than $1$, uniformly in that range. Hence this part of
the sum is bounded by a constant times
$\gamma^M/\binom nM=O(n^{-M})$. For $n/4\leq j\leq n/2$, use
\[
 \binom nj
 =\prod_{k=0}^{j-1}\frac{n-k}{j-k}
 \geq\left(\frac nj\right)^j\geq2^j.
\]
The sum of $\gamma^j/\binom nj$ over this second range is therefore
exponentially small, since $\gamma/2<1$. Symmetry in $j$ and
$n-j$ proves the binary convolution estimate.

Taking $M=0$ shows that the unrestricted binary convolution is
bounded by a constant times $u_n$. Increasing that constant to
cover the finitely many small values of $n$, induction gives
\[
 \sum_{n_1+\cdots+n_q=n}u_{n_1}\cdots u_{n_q}
 \leq C_q u_n
\]
for every fixed $q$. In a restricted $q$-tuple, choose an index
whose value $j$ is smallest. Then $M\leq j\leq n/q\leq n/2$.
There are at most $q$ choices of that index. After fixing it,
bound the convolution of the remaining $q-1$ indices by the
unrestricted estimate. The resulting upper bound is
\[
 qC_{q-1}\sum_{j=M}^{\lfloor n/2\rfloor}u_j u_{n-j}
 =O(u_n n^{-M}),
\]
as required. For $q=2$, the factor $C_{q-1}$ can be taken to be
$1$.
\end{proof}

The coefficient linearization in the next theorem is related to
Borinsky's asymptotic linearization~\cite[Proposition~55 and
Remark~57]{Borinsky18}. The
argument given here requires only a positive coefficient ratio
limit; it does not assume a complete asymptotic expansion.

\begin{theorem}
\label{thm:ratio-algebraic}
Let $f(t)=\sum_{n\geq0}c_n t^n\in\Q[[t]]$, where $c_n>0$ for
all sufficiently large $n$. Suppose that
\[
 \lim_{n\to\infty}\frac{c_{n+1}}{nc_n}=K,
 \qquad 0<K<\infty.
\]
If $f$ is differentially algebraic over $\C(t)$, then $K$ is an
algebraic number.
\end{theorem}

\begin{proof}
The proof has three steps. First, we reduce to a differential
polynomial over \(\Q[t]\) whose total degree \(d\) in the
\(Y\) variables is minimal. Second, we expand this relation in a
Taylor expansion around the truncated series \(f_N\). The linear
part carries the dominant contribution \(c_n n^E\), while the
nonlinear remainder is smaller by \cref{lem:factorial-convolution}. Third, comparing
coefficients gives a nonzero polynomial equation for \(K\).

We first show that a differential polynomial relation can be
chosen over $\Q[t]$. Clear the denominators in a relation over
$\C(t)$ and retain the finitely many monomials in
$t,Y_0,\ldots,Y_r$ that occur. Regard their coefficients as
unknowns. Substitution of $f,f',\ldots,f^{(r)}$ and comparison of
coefficients of $t$ gives a homogeneous linear system with
rational entries and finitely many columns. Its row space has a
basis consisting of finitely many rows. Consequently its kernel
is already the kernel of a finite rational matrix, whose rank is
the same over $\Q$ and $\C$. A nonzero complex solution therefore
implies a nonzero rational solution.

Among all nonzero polynomials
$P\in\Q[t,Y_0,\ldots,Y_r]$ satisfying
\[
 P(t,f,f',\ldots,f^{(r)})=0,
\]
choose one whose total degree in the $Y$ variables is smallest.
Denote that degree by $d$. It is positive, since a nonzero
polynomial in $t$ alone cannot vanish as a formal series. For
every index $i$ on which $P$ depends, its partial derivative with
respect to $Y_i$ is a nonzero polynomial of smaller total degree.
Minimality of $d$ therefore implies that
\[
 H_i(t)=\frac{\partial P}{\partial Y_i}
           (t,f,f',\ldots,f^{(r)})
\]
is a nonzero formal series. Write its first term as
\[
 H_i(t)=h_i t^{\nu_i}+O(t^{\nu_i+1}),
 \qquad h_i\in\Q\setminus\{0\},\quad \nu_i\geq0,
\]
and set $E=\max_i(2i-\nu_i)$, where the maximum is over the
indices on which $P$ depends. The choice of \(E=\max_i(2i-\nu_i)\) records the dominant power of
\(n\) in the linear term below. The truncation index \(M\) will be
chosen large enough to make the nonlinear remainder negligible
compared with \(c_n n^E\).

The ratio hypothesis implies, by a finite product of consecutive ratios, that for every fixed integer $s$,
\[
 \frac{c_{n+s}}{c_n}\sim K^s n^s.
\]
Negative shifts follow by taking reciprocals of the corresponding
positive-shift products. In particular, for fixed nonnegative
integers $i,k$,
\[
 \frac{[t^{n-k}]f^{(i)}}{c_n}
 =\frac{(n-k+i)!}{(n-k)!}\frac{c_{n-k+i}}{c_n}
 \sim K^{i-k}n^{2i-k}.
\]

Choose a large fixed integer $N$, and set
\[
 f_N(t)=\sum_{j=0}^{N-1}c_j t^j,
 \qquad R(t)=f(t)-f_N(t),\qquad M=N-r.
\]
Taylor expansion of the polynomial $P$ in its $Y$ variables is an
exact finite identity:
\begin{align*}
 0={}&P(t,f_N,f_N',\ldots,f_N^{(r)})
       +\sum_i H_{i,N}(t)R^{(i)}
       +Q_N(t,R,R',\ldots,R^{(r)}),\\
 H_{i,N}(t)={}&\frac{\partial P}{\partial Y_i}
                    (t,f_N,f_N',\ldots,f_N^{(r)}).
\end{align*}
Here each $H_{i,N}$ is a polynomial, and every monomial of $Q_N$
has at least two factors among $R,R',\ldots,R^{(r)}$. The
coefficients of those monomials are polynomials in $t$.

Every difference $f^{(i)}-f_N^{(i)}$ is divisible by $t^M$ for
$0\leq i\leq r$. Substitution in a polynomial thus shows that
$H_i-H_{i,N}$ is divisible by $t^M$. Choose $N$ so large that
\[
 M>\max_i\nu_i,\qquad M>2rd-E,
\]
and so that all coefficients $c_j$ with $j\geq M$ are positive.
It follows that $H_{i,N}$ has first term $h_i t^{\nu_i}$.
If $P$ does not depend on $Y_i$, then both $H_i$ and $H_{i,N}$
are identically zero.

For sufficiently large $n$, the coefficient of
$P(t,f_N,f_N',\ldots,f_N^{(r)})$ is zero, since it is a
polynomial. Likewise, in the linear part, the relevant
coefficients of $R^{(i)}$ equal those of $f^{(i)}$. Each
$H_{i,N}$ contains only finitely many terms. Applying the
fixed-shift estimate to these terms gives
\[
 [t^n]\sum_i H_{i,N}R^{(i)}
 =c_n n^E
   \left(\sum_{2i-\nu_i=E}h_iK^{i-\nu_i}+o(1)\right).
\]
Indeed, the first term of $H_{i,N}$ contributes
$h_iK^{i-\nu_i}c_n n^{2i-\nu_i}(1+o(1))$, and every later
term has a smaller power of $n$.

We now bound the nonlinear remainder. Replace any nonpositive
initial coefficients of $(c_n)$ by positive numbers, obtaining a
sequence $(u_n)$ that agrees with $(c_n)$ at all sufficiently
large indices. It has the same ratio limit, so \cref{lem:factorial-convolution} applies. Enlarge $N$ if
necessary so that $u_j=c_j$ for $j\geq M$. The fixed-shift
estimate also gives
\[
 \bigl|[t^j]R^{(i)}\bigr|
 =\frac{(j+i)!}{j!}c_{j+i}
 \leq C_i u_j j^{2i}
\]
whenever this coefficient is nonzero. All such indices satisfy
$j\geq N-i\geq M$, so the constants can be chosen uniformly
over the relevant indices.

Consider a remainder monomial containing $q\geq2$ factors, with
derivative orders $i_1,\ldots,i_q$. Its coefficient of degree $n$
is bounded in absolute value by a constant times
\[
 n^{2(i_1+\cdots+i_q)}
 \sum_{\substack{n_1+\cdots+n_q=n\\n_1,\ldots,n_q\geq M}}
 u_{n_1}\cdots u_{n_q}
 =O(c_n n^{2(i_1+\cdots+i_q)-M}).
\]
There are at most $d$ factors, each of derivative order at most
$r$, so this is $O(c_n n^{2rd-M})$. Multiplication by a fixed
power $t^p$, with $p\geq0$, preserves that bound: the coefficient
is shifted from $n$ to $n-p$, and
$u_{n-p}/c_n=O(n^{-p})$. Since the remainder contains only
finitely many such terms,
\[
 [t^n]Q_N(t,R,R',\ldots,R^{(r)})
 =O(c_n n^{2rd-M})=o(c_n n^E).
\]

The exact Taylor identity therefore gives
\[
 \sum_{2i-\nu_i=E}h_iK^{i-\nu_i}=0.
\]
This is a nonzero Laurent polynomial equation over $\Q$.
To check that its polynomial is nonzero, observe that
$\nu_i=2i-E$ in this sum, so the exponent $i-\nu_i=E-i$
is different for each index $i$. Every coefficient $h_i$ is
nonzero. Multiplying by a sufficiently large power of $K$
produces a nonzero ordinary polynomial over $\Q$ satisfied by
$K$. Thus $K$ is algebraic.
\end{proof}

\subsection{Proof of \cref{thm:block-egfs} by ratio limits}
\label{subsec:ratio-limits}

The criterion \cref{thm:ratio-algebraic} requires the limit of the
consecutive ratio $c_{n+1}/(nc_n)$. Our growth estimates in
\cref{thm:ABnotDfinite} give only the $n$-th root limits of the
normalized counts. We therefore need a passage from root limits to
ratio limits. The passage uses the following properties of $d_n$ and
$M_n=e_{n+1}/((2n+2)(2n+1))$.

\begin{lemma}
\label{lem:moment-facts}
The sequences \(d_n\) and \(M_n\) are moment sequences of positive
finite measures supported on \([0,L^2]\), and
\(e_{n+1}=(2n+2)(2n+1)M_n\). Moreover, both sequences are
log-convex:
\[
 d_n^2\le d_{n-1}d_{n+1},\qquad
 M_n^2\le M_{n-1}M_{n+1}
 \qquad(n\ge1).
\]
\end{lemma}

The proof of \cref{lem:moment-facts} is given in
\cref{subsec:moments}. It uses a probabilistic representation of
$d_n$ and $M_n$ as moments of positive measures, which is
independent of the coefficient criterion above. 

Now we are ready to give the proof of \cref{thm:block-egfs} as follows. 

Log-convexity implies that the ratios \(d_{n+1}/d_n\) and
\(M_{n+1}/M_n\) increase with \(n\). Since the representing
measures are supported on \([0,L^2]\), we have \(d_n\le L^{2n}\) and
\(M_n\le L^{2n}\), so these ratios are bounded above by \(L^2\).
Hence both limits exist. For a positive sequence \((m_n)\) with
\(m_n=m_0\prod_{j=0}^{n-1}(m_{j+1}/m_j)\), the root limit
\(m_n^{1/n}\) equals the ratio limit whenever the latter exists.
The root limits are known by \cref{thm:ABnotDfinite}: \(d_n^{1/n}\to L^2\), and
\(M_n^{1/n}\to L^2\) follows from
\(e_{n+1}=(2n+2)(2n+1)M_n\) and \(e_n^{1/n}\to L^2\). Therefore
\[
 \frac{d_{n+1}}{d_n}\longrightarrow L^2,
 \qquad
 \frac{M_n}{M_{n-1}}\longrightarrow L^2.
\]
Consequently
\[
 \frac{e_{n+1}}{e_n}
 =\frac{(2n+2)(2n+1)}{(2n)(2n-1)}
      \frac{M_n}{M_{n-1}}
 \longrightarrow L^2.
\]

For the coefficients $c_n=a_n/n!$ of $F_2$, the factorial
normalization gives the exact ratio
\[
 \frac{c_{n+1}}{nc_n}
 =\frac{d_{n+1}}{d_n}
  \frac{(4n+1)(4n+2)(4n+3)(4n+4)}
       {(2n+1)(2n+2)n(n+1)}
 \longrightarrow64L^2=2\pi^2.
\]
For $c_n=b_n/n!$, the same calculation uses $e_n$ in place
of $d_n$ and gives the same limit. Both coefficient sequences
consist of positive rational numbers. Since $2\pi^2$ is
transcendental by the transcendence of $\pi$ (see \cite{Lindemann82}),
\cref{thm:ratio-algebraic} proves that
$F_2$ and $F_3$ are differentially transcendental. Every
$D$-finite series satisfies a differential polynomial equation,
so neither series is $D$-finite.

It remains to verify the stated convergence properties.
For either coefficient sequence, the ratio just proved implies
\[
 \frac{c_{n+1}/(n+1)!}{c_n/n!}
 =\frac{c_{n+1}}{(n+1)c_n}\longrightarrow2\pi^2.
\]
Taking a finite product of these ratios and then taking the $n$-th roots gives
\[
 \left(\frac{c_n}{n!}\right)^{1/n}\longrightarrow2\pi^2.
\]
The Cauchy--Hadamard formula therefore gives radius
$1/(2\pi^2)$ for each ordinary Borel transform. For any
$A>2\pi^2$, the same limit gives $c_n\leq CA^n n!$ after
increasing $C$ to cover the initial indices. Finally,
\[
 c_n^{1/n}
 =\left(\frac{c_n}{n!}\right)^{1/n}(n!)^{1/n}
 \longrightarrow\infty,
\]
because $(n!)^{1/n}\to\infty$. Thus the radii of convergence
of $F_2$ and $F_3$ are zero. This completes the proof of
\cref{thm:block-egfs} modulo the proof of \cref{lem:moment-facts} in \cref{subsec:moments}.

\subsection{Positive moment representations}
\label{subsec:moments}

Recall the normalized counts and their exponential growth limits:
\[
 d_n=\frac{(2n)!}{(4n)!}a_n,
 \qquad e_n=\frac{(2n)!}{(4n)!}b_n,
 \qquad
 d_n^{1/n},\ e_n^{1/n}\longrightarrow L^2,
 \qquad L=\frac{\pi}{4\sqrt2}.
\]
The limits were proved in \cref{thm:ABnotDfinite}.
We strengthen them to ratio limits by representing the normalized
counts as moments of positive measures.

The construction below applies to the equations for \(S\) and \(V\).
It uses three features that are specific to \(q(x)=(1-x^2)/2\) and
to the initial-value problems \eqref{eq:S-ivp} and \eqref{eq:V-ivp}.
First, for \(S\), after the substitution \(z=w^2\), the functions
\(U=S_x\) and \(W=wvS\) satisfy a \(2\times2\) linear system whose
first equation has no zeroth-order term. For \(V\), the same system
has an additional forcing term \(w^2(1-x)\) in the first equation.
Second, the coefficient
\[
 \kappa(x)=-\frac{q'(x)}{4q(x)}=\frac{x}{2(1-x^2)}
\]
is nonnegative on \([0,1)\), so it can serve as a switching rate.
Third, the integral \(\int_0^x v(u)\,du\) is bounded by \(L\), so the
resulting random variable has bounded support. Under these three
conditions, the \(2\times2\) system is the generator of a two-state
Markov process, and its solution is the exponential moment of that
process. No general differential equation admits such a
representation.

Set
\[
 v(x)=\sqrt{q(x)},\qquad
 \kappa(x)=\frac{x}{2(1-x^2)},\qquad 0\leq x<1.
\]
In particular, $v'/v=-2\kappa$, and $\kappa$ is nonnegative
and continuous on every compact subinterval of $[0,1)$.
We give the switching-process construction and the elementary
differential equations for its weighted probabilities explicitly.

Let $\Pi(u)$, $u\geq0$, be a Poisson process of rate $1$,
starting at $\Pi(0)=0$ and having right-continuous counting paths: its
increments on disjoint intervals are independent, and
\[
 \Pr\{\Pi(u+h)-\Pi(u)=j\}=\ee^{-h}\frac{h^j}{j!}
 \qquad(j\geq0).
\]
Let $\varepsilon$ be independent of this process, with values
$1$ and $-1$, each of probability $1/2$. Define the common
time-changed process
\[
 N(x)=\Pi\left(\int_0^x\kappa(u)\dd u\right),
 \qquad 0\leq x<1.
\]
For $0\leq s\leq x<1$, set
\[
 \sigma_{s,x}
 =\varepsilon(-1)^{N(x)-N(s)},
 \qquad
 X_{s,x}=\int_s^x v(u)\sigma_{s,u}\dd u.
\]
This common construction makes the random integrals jointly
measurable in $s,x$, as required for the double integrals below.
On each compact interval $[s,c]$ with $c<1$, the cumulative
Poisson time is finite. Thus $\sigma_{s,x}$ switches sign a
finite number of times almost surely on that interval, and its
switching rate at $x$ is $\kappa(x)$. The time change is explicit:
\[
 \int_s^x\kappa(u)\dd u
 =\frac14\log\frac{1-s^2}{1-x^2}.
\]
Although this quantity tends to infinity as $x\to1$, the integral
defining $X_{s,x}$ has a pathwise limit $X_{s,1}$, because
\[
 |X_{s,x'}-X_{s,x}|
 \leq\int_x^{x'}v(u)\dd u
 \qquad(s\leq x\leq x'\leq1).
\]
Moreover,
\[
 |X_{s,x}|\leq\int_s^xv(u)\dd u\leq L.
\]
No value of the terminal sign $\sigma_{s,1}$ is needed.

For $w\in\C$ and $s\leq x<1$, define
\[
 p_\pm(s,x;w)
 =\mathbb E\left[
       \ee^{wX_{s,x}}\mathbf1_{\{\sigma_{s,x}=\pm1\}}
                 \right],
\]
where $\mathbb E$ denotes expectation and $\mathbf1$ an indicator.
Their initial values are $p_\pm(s,s;w)=1/2$. We claim that
\begin{equation}
 \begin{split}
 \partial_xp_+&=(wv-\kappa)p_++\kappa p_-,\\
 \partial_xp_-&=\kappa p_++(-wv-\kappa)p_-.
 \end{split}
 \label{eq:weighted-switching}
\end{equation}
To verify this, fix $x$ in a compact interval below $1$ and
condition on the sign at $x$ and the Poisson increments over
$[x,x+h]$. The probabilities of no switch, one switch, and at
least two switches are, respectively,
\[
 1-\kappa(x)h+O(h^2),\qquad
 \kappa(x)h+O(h^2),\qquad O(h^2).
\]
On a path with no switch, the extra exponential weight is
$1\pm wv(x)h+O(h^2)$, according to the sign. On a path with
one switch it is $1+O(h)$; its $O(h)$ correction is multiplied
by a switching probability of order $h$. The weights already
accumulated before $x$ are bounded by $\ee^{|w|L}$. Independent
Poisson increments therefore give
\begin{align*}
 p_+(s,x+h;w)
 &=p_+(s,x;w)
   +h\bigl((wv-\kappa)p_++\kappa p_-\bigr)+O(h^2),\\
 p_-(s,x+h;w)
 &=p_-(s,x;w)
   +h\bigl(\kappa p_++(-wv-\kappa)p_-\bigr)+O(h^2),
\end{align*}
with the coefficients on the right evaluated at $x$.
Division by $h$ and passage to the limit proves
\eqref{eq:weighted-switching}.

\begin{proposition}
\label{prop:block-moments}
The first normalized count is the moment
\[
 d_n=\mathbb E X_{0,1}^{2n}\qquad(n\geq0).
\]
For the second family, set
\[
 M_j=\int_{0<s<x<1}(1-s)\mathbb E X_{s,x}^{2j}
                         \dd s\dd x\qquad(j\geq0).
\]
Then
\[
 e_n=(2n)(2n-1)M_{n-1}\qquad(n\geq1).
\]
Both $(d_n)$ and $(M_n)$ are moment sequences of positive
finite measures supported on $[0,L^2]$.
\end{proposition}

\begin{proof}
For the first family, evaluate \eqref{eq:S-ivp} at $z=w^2$ and
set $U=S_x$ and $W=wvS$, with dependence on $x,w$ understood.
The differential equation for $S$ and the identity
$v'/v=-2\kappa$ give
\[
 U'=wvW,\qquad W'=wvU-2\kappa W.
\]
Consequently $U+W$ and $U-W$ satisfy the same system as
$2p_+(0,x;w)$ and $2p_-(0,x;w)$ in
\eqref{eq:weighted-switching}. Their initial values are also
the same, since $U(0)=1$ and $W(0)=0$. Uniqueness for the
linear system on each compact interval below $1$ gives
\[
 S_x(x,w^2)=U(x)=p_+(0,x;w)+p_-(0,x;w)
                    =\mathbb E\ee^{wX_{0,x}}.
\]
The solution $S_x$ is continuous at $x=1$ by
\eqref{eq:S-ivp}. The random integral converges pathwise there,
and its exponential is bounded by $\ee^{|w|L}$. Dominated
convergence and \eqref{eq:boundary-values} therefore yield
\[
 A(w^2)=\mathbb E\ee^{wX_{0,1}}.
\]
Changing the independent initial sign $\varepsilon$ to
$-\varepsilon$ changes $X_{0,1}$ to $-X_{0,1}$ without changing
its distribution. Thus all odd moments vanish. The exponential
series converges uniformly for $|X_{0,1}|\leq L$ on bounded
sets of $w$, so expectation can be taken term by term:
\[
 A(w^2)=\sum_{n\geq0}
             \frac{\mathbb E X_{0,1}^{2n}}{(2n)!}w^{2n}.
\]
Comparison with the definition of $A$ proves the asserted
formula for $d_n$.

For the second family, evaluate \eqref{eq:V-ivp} at $z=w^2$
and set $U=V_x$ and $W=wvV$. This time the equations are
\[
 U'=wvW+w^2(1-x),\qquad W'=wvU-2\kappa W.
\]
Hence $U+W$ and $U-W$ satisfy
\eqref{eq:weighted-switching} with the same extra term
$w^2(1-x)$ in both equations, and with zero initial values.
It follows that
\[
 U(x)\pm W(x)
 =2w^2\int_0^x(1-s)p_\pm(s,x;w)\dd s.
\]
This identity can be checked directly: differentiation of the
right side gives the homogeneous system under the integral,
and its upper endpoint contributes
$2w^2(1-x)p_\pm(x,x;w)=w^2(1-x)$; its value at $x=0$ is
zero. Linear-system uniqueness then proves the identity.
Adding the two signs and dividing by $2$ gives
\[
 V_x(x,w^2)
 =w^2\int_0^x(1-s)\mathbb E\ee^{wX_{s,x}}\dd s.
\]
Integrate over $0<x<1$ and use
\eqref{eq:boundary-values}:
\[
 B(w^2)-1
 =w^2\int_{0<s<x<1}(1-s)\mathbb E\ee^{wX_{s,x}}
                                      \dd s\dd x.
\]
All weights are bounded by $\ee^{|w|L}$, and the integration
domain has finite area. The same symmetry and uniform
convergence argument as before therefore gives
\[
 B(w^2)-1
 =\sum_{j\geq0}\frac{M_j}{(2j)!}w^{2j+2}.
\]
For $n\geq1$, comparison of coefficients yields
\[
 \frac{b_n}{(4n)!}=\frac{M_{n-1}}{(2n-2)!},
\]
which is equivalent to the asserted formula for $e_n$.
In particular, $M_0=\int_0^1(1-s)^2\dd s=1/3$,
in agreement with $b_1/4!=1/3$.

Finally, $d_n$ is the $n$-th moment of the distribution of
$X_{0,1}^2$. The sequence $M_n$ is the $n$-th moment of the
positive measure obtained by integrating the distributions of
$X_{s,x}^2$ with weight $(1-s)\dd s\dd x$. Their total masses
are $1$ and $1/3$, respectively. The bound $|X_{s,x}|\leq L$
places both supports in $[0,L^2]$.
\end{proof}

We also record the log‑convexity used in \cref{subsec:ratio-limits}.
For a positive finite measure supported on \([0,L^2]\), the
Cauchy--Schwarz inequality gives
\[
 m_n^2
 =\left(\int x^n\dd\mu\right)^2
 \le\left(\int x^{n-1}\dd\mu\right)
     \left(\int x^{n+1}\dd\mu\right)
 =m_{n-1}m_{n+1}.
\]
Applied to the measures representing \(d_n\) and \(M_n\), this gives
\(d_n^2\le d_{n-1}d_{n+1}\) and \(M_n^2\le M_{n-1}M_{n+1}\) for
\(n\ge1\).

\appendix
\section{Subset sum formulas for \texorpdfstring{$a_n$}{an} and \texorpdfstring{$b_n$}{bn}}
\label{sec:subset-sums}

The elementary simplex integral used here is
\begin{equation}
 \int_{0<u_1<\cdots<u_m<1}
       u_1^{r_1}\cdots u_m^{r_m}\dd u_1\cdots\dd u_m
 =\prod_{k=1}^m\frac{1}{k+r_1+\cdots+r_k},
 \qquad r_k\geq0.
 \label{eq:monomial-simplex}
\end{equation}
To verify it, integrate $u_1$ first. The factor
$(r_1+1)^{-1}$ is obtained, and the exponent of $u_2$ increases by
$r_1+1$. The next integration gives the factor
$(r_1+r_2+2)^{-1}$. Repeating this operation through $u_m$ gives
exactly the product in \eqref{eq:monomial-simplex}.

\begin{proposition}\label{prop:subset-sums}
For every integer $n\geq0$,
\begin{equation}
 a_n=\frac{(4n)!}{2^n}
 \sum_{\varepsilon\in\{0,1\}^n}
 \frac{(-1)^{\varepsilon_1+\cdots+\varepsilon_n}}
 {\displaystyle\prod_{j=1}^n
   \left(2j-1+2\sum_{k=1}^{j-1}\varepsilon_k\right)
   \left(2j+2\sum_{k=1}^{j}\varepsilon_k\right)}.
 \label{eq:a-sum}
\end{equation}
For every integer $n\geq1$,
\begin{equation}
 b_n=\frac{(4n)!}{2^{n-1}}
 \sum_{\varepsilon\in\{0,1\}^n}
 \frac{(-1)^{\varepsilon_1+\cdots+\varepsilon_n}}
 {\displaystyle\prod_{j=1}^n
   \left(2j-1+\varepsilon_1+2\sum_{k=2}^{j}\varepsilon_k\right)
   \left(2j+\varepsilon_1+2\sum_{k=2}^{j}\varepsilon_k\right)}.
 \label{eq:b-sum}
\end{equation}
Empty products are $1$, and sums over empty index sets are $0$.
The value at $n=0$ in \eqref{eq:a-sum} is interpreted using the single
empty binary vector.
\end{proposition}

\begin{proof}
Expand the integrand in \eqref{eq:a-density} as
\[
 2^{-n}\prod_{j=1}^n(1-x_{2j}^2)
 =2^{-n}\sum_{\varepsilon\in\{0,1\}^n}
    (-1)^{\varepsilon_1+\cdots+\varepsilon_n}
    \prod_{j=1}^n x_{2j}^{2\varepsilon_j}.
\]
In \eqref{eq:monomial-simplex}, the cumulative exponent at position
$2j-1$ is $2\sum_{k<j}\varepsilon_k$, and the cumulative exponent at
position $2j$ is $2\sum_{k\leq j}\varepsilon_k$. Substituting these
values proves \eqref{eq:a-sum}.

For \eqref{eq:b-density}, the expansion is
\[
 2^{-(n-1)}(1-x_1)\prod_{j=2}^n(1-x_{2j-1}^2) =2^{-(n-1)}\sum_{\varepsilon\in\{0,1\}^n}(-1)^{\varepsilon_1+\cdots+\varepsilon_n}
 x_1^{\varepsilon_1}\prod_{j=2}^n x_{2j-1}^{2\varepsilon_j}.
\]
At positions $2j-1$ and $2j$, the cumulative exponent is the same,
namely $\varepsilon_1+2\sum_{k=2}^j\varepsilon_k$.
Equation~\eqref{eq:monomial-simplex} now gives \eqref{eq:b-sum}.
\end{proof}




\noindent\textbf{Data availability:} Data availability is not applicable to this article as no new data were created or analyzed in this study.


\vspace{5mm}
\noindent
{\small \textbf{Acknowledgments:}}
We are also grateful to Dr. Zihao Zhang for many useful suggestions.
Feihu Liu was partially supported by the Postdoctoral Fellowship Program and China Postdoctoral Science Foundation (No. BX2026002).
Guoce Xin was partially supported by the National Natural Science Foundation of China (No. 12571355).


\end{document}